\documentclass[11pt]{article}
\usepackage[utf8]{inputenc}
\usepackage{caption}
\usepackage[T1]{fontenc}
\usepackage{graphicx}
\usepackage[english]{babel}
\usepackage{csquotes}
\usepackage{amssymb}
\usepackage{amsthm}
\usepackage{amsbsy}
\usepackage[version=4]{mhchem} %notacion quimica /ce
\usepackage{hyperref}
\usepackage{placeins} %para el floatbarrier
\usepackage{url}
\usepackage{multirow}
\usepackage{amsmath}
\usepackage{tikz}
\usepackage[sc]{mathpazo}
\usepackage[letterpaper,top=2.5cm,bottom=2cm,left=1.3cm,right=1.5cm,marginparwidth=1.75cm]{geometry}
\usepackage[colorinlistoftodos]{todonotes}
\usepackage{units}
\usepackage[shortlabels]{enumitem}
\usepackage{float}
\usepackage{blindtext}
\usepackage{microtype}
\usepackage{makecell}
\usepackage{lettrine}
\usepackage{titling}
\usepackage{booktabs}
\usepackage{ragged2e}
\usepackage{cuted}
\usepackage{authblk}
\usepackage{adjustbox}
\usepackage{xcolor}
\usepackage{lipsum}
\usepackage[normalem]{ulem}
\usepackage{lscape}
\usepackage{tabularx}
\usepackage{cancel}
\usepackage{lmodern}%get scalable font
\usepackage{dsfont}
\usepackage{lipsum}
\usepackage{relsize}
\usepackage{comment}
\usepackage{physics} % para la notation bra y ket
\usepackage{mathtools} % para hacer la flecha en ambas direcciones
\usepackage{tablefootnote}
\usepackage{derivative}%ecuaciones diferenciales de grado n usar \odv[n]{f}{x} y parciales de grado n \pdv[n]{f}{x}
\setlist[itemize]{noitemsep}
\usepackage{abstract}

\providecommand{\keywords}[1]
{
  \small	
  \textbf{\textit{Keywords---}} #1
}

\usepackage{titlesec}% Allows customization of titles

\titleformat{\section}[hang]{\large\centering}{\thesection.}{1em}{} % Change the look of the section titles

\titleformat{\subsubsection}[block]{\large}{\thesubsubsection.}{1em}{} % Change the look of the subsubsection titles

\usepackage[labelfont=bf,textfont=it]{caption}
\usepackage{subcaption}

\titleformat{\subsection}[block]{\large}{\thesubsection.}{1em}{}
\usepackage{hyperref}
\usepackage{listings}
\hypersetup{
    colorlinks=true,
    citecolor=blue,
    filecolor=blue,      
    urlcolor=blue,
    allcolors=blue
}

\usepackage{cleveref}
\usepackage[
	backend = biber,
	bibencoding = utf8,
	sorting = none
]{biblatex}

\makeatletter 	  	 	
\renewcommand*{\@fnsymbol}[1]{\ensuremath{\ifcase#1\or *\or \dagger\or \ddagger\or
    \mathsection\or \mathparagraph\or \|\or **\or \dagger\dagger
    \or \ddagger\ddagger \else\@ctrerr\fi}}
\makeatother

\newcommand{\subtitle}[1]{%
  \posttitle{%
    \par\end{center}
    \begin{center}\large#1\end{center}
    \vskip0.5em}%
}

\usepackage{appendix}
\addto\captionsspanish{%

}
\title{Higher order derivatives of the adjugate matrix and the Jordan form}
\author[1]{Jorge I. Rubiano-Murcia}
\affil[1]{Departamento de F\'isica, Universidad Nacional de Colombia, Carrera 45 No. 26-85, Edificio Uriel Guti\'errez, Bogot\'a D.C., Colombia. Email {\tt jrubianom@unal.edu.co}}
\author[2]{Juan Galvis}
\affil[2]{Departamento de Matem\'aticas, Universidad Nacional de Colombia, Carrera 45 No. 26-85, Edificio Uriel Guti\'errez, Bogot\'a D.C., Colombia. Email {\tt jcgalvis@unal.edu.co}}

\newtheorem{theorem}{Theorem}[section]
\newtheorem{corollary}{Corollary}[theorem]
\newtheorem{lemma}[theorem]{Lemma}
\theoremstyle{definition}

\date{\today}
\begin{document}
\maketitle
\begin{abstract}
In this short note, we show that the higher-order derivatives of the adjugate matrix $\mbox{Adj}(z-A)$, are related to the nilpotent matrices and projections in the Jordan decomposition of the matrix $A$. These relations appear as a factorization of the derivative of the adjugate matrix as a product of factors related to the eigenvalues, nilpotent matrices and projectors.  
The novel relations are obtained using the Riesz projector and functional calculus. 
The results presented here can be considered to be a generalization of Thompson and McEnteggert's theorem relating the adjugate matrix to the orthogonal projection on the eigenspace of simple eigenvalues for symmetric matrices. They can also be seen as a complement to some earlier results by B. Parisse, M. Vaughan that relate derivatives of the adjugate matrix to the invariant subspaces associated with an eigenvalue. Our results can also be interpreted as a general eigenvector-eigenvalue identity.
Many previous works have dealt with relations between the projectors on the eigenspaces and the derivatives of the adjugate matrix with the characteristic spaces but it seems that there is no explicit mention in the literature of the factorization of the higher-order derivatives of the adjugate matrix as a matrix multiplication involving nilpotent and projector matrices, which appear in the Jordan decomposition theorem.
\end{abstract}\hspace{10pt}


%\tableofcontents
\section{Introduction and main results}
Let $A \in \mathbb{C}^{n\times n}$ be a square complex matrix of order $n$ and denote its adjugate matrix by $\mbox{Adj}(A)$. We recall a result of  {Thompson and McEnteggert}  (\cite{iserles2002acta,parlett1998symmetric,thompson1968principal, Castillo2021OnAF}) which states that for a Hermitian matrix $A$ whose characteristic polynomial is $p(z)$, if $\lambda_i$ is a simple eigenvalue of $A$ and $z_i$ is a corresponding unit eigenvector, then
\begin{equation}
    \mbox{\normalfont Adj}(\lambda_i - A) = \frac{d p(z)}{dz}\Bigr|_{\lambda_i} z_i z_i^*\,.    
    \label{eq:def:TM}
\end{equation}
This means that $\mbox{\normalfont Adj}(\lambda_i - A)$ is a scalar multiple of the orthogonal projection on the characteristic eigenspace related to $\lambda_i$. Denote by $n_i$ the (algebraic) multiplicity of the eigenvalue $\lambda_i$, where $1\leq i\leq m$ and $m$ is the total number of distinct eigenvalues. 
It is also known that the matrix
\[
    \frac{d^{n_i-1}  \mbox{\normalfont Adj}(z - A)}{dz^{n_i-1}}\Bigr|_{\lambda_i}
\]
spans the root (characteristic) subspace associated with $\lambda_i$. 
See \cite{1parisse:hal-00003444} and the references therein. In this short article, we extend and unify these two observations. In particular, we extend 
\eqref{eq:def:TM} to the case of repeated eigenvalues, higher derivatives and general matrices. 
We also show that the higher-order derivatives of the adjugate matrix $\mbox{Adj}(z-A)$, are related to the nilpotent matrices and projections in the Jordan decomposition of the matrix $A$.
The following factorization is obtained in Section \ref{sec:derivatives}, 
\begin{equation}\label{eq:spoiler}
    \frac{d^{n_i-1}  \mbox{\normalfont Adj}(z - A)}{dz^{n_i-1}}\Bigr|_{\lambda_i} = (n_i-1)! \, \textbf{   } 
    \prod_{j\not =i }^m (N_i+\lambda_i -\lambda_j)^{n_j} P_i\,,
\end{equation}
where $N_i$ and $P_i$ are matrices present in Jordan form. Specifically, the matrix $N_i$ is an $n\times n$ nilpotent matrix of degree less than or equal to $n_i$ ($N_i^{n_i}=0$). On the other hand, $P_i$ is the oblique projection on the root subspace associated with $\lambda_i$. It satisfies $P_i^2=P$ and the columns of $P_i$ span the root subspace. For Hermitian matrices $N_i=0$, thus, we recover the result of Thompson and McEnteggert. In the general case, we can write
\[
   N_i=N_i P_i = \frac{1}{(n_i-2)!}\frac{d^{n_i-2} }{dz^{n_i-2}} \left( \frac{\mbox{\normalfont Adj}(z - A)}{\prod_{j\not =i }^m (z-\lambda_j)^{n_j}}\right)\Bigr|_{\lambda_i}\,,
\]
which is a nilpotent matrix of degree less than or equal to $n_i$. The matrix $N_i$ is related to the action of $A$ restricted to the root subspace associated with 
the eigenvalue $\lambda_i$. Furthermore, we discuss additional related results. Among them, we prove that for all $i$ and $0\leq s< n_i-1$, we have
    \[
    \frac{d^{s} \mbox{\normalfont Adj}(z - A)}{dz^{s}}\Bigr|_{\lambda_i} = s! N_i^{n_i-1-s}\, \textbf{   } \prod_{j\not =i }^m (N_i+\lambda_i-\lambda_j)^{n_j} P_i.
    \]
This formula is consistent with \eqref{eq:spoiler} when 
$s=n_i-1$ and $N_i\not=0$. Note that the above novel identities relate the eigenvalues of $A$, $\{\lambda_j\}_{j=1,\dots , m}$, the eigenvectors - as $P_i$ spans the root spaces associated with $\lambda_i$ - and the derivatives of the cofactors evaluated at the eigenvalue $\lambda_i$. As a result, these identities can be interpreted as generalized \emph{eigenvector-eigenvalue identities}. For an overview of eigenvector-eigenvalue identities, see the survey in \cite{denton2022eigenvectors}. It is observed that \cite{denton2022eigenvectors} derives eigenvector-eigenvalue identities via various methods, including the formula \eqref{eq:def:TM} by Thompson and McEnteggert.

The following sections of the paper are organized as follows. Section \ref{sec:prelim} presents background material and related previous work. Section \ref{sec:derivatives} contains a detailed presentation of our main results. Section \ref{sec:examples} provides some illustrative examples, and Section \ref{sec:conclusions} concludes this document with final remarks.

\section{Preliminaries and related works}\label{sec:prelim}
Let $A \in \mathbb{C}^{n\times n}$ be a complex matrix and let $\Gamma \subset \mathbb{C}$ be a positively oriented rectifiable Jordan curve (\cite[chapter 3,pp 68]{ahlfors1953complex}) that does not contain eigenvalues of $A$ and encloses a region $D$. Then, the \textit{Riesz projector} is defined as, see \cite[chapter 1,pp 3-6]{gohberg1978introduction}
\begin{align}
    P_\Gamma &= \frac{1}{2\pi i} \oint\limits_{\Gamma} (z-A)^{-1} \, \mathrm{d}z\, .
    \label{eq:RieszProjector}
\end{align}
When $D$ contains only one eigenvalue of $A$, say $\lambda_i$, the Riesz projector is denoted by $P_{\lambda_i}$ or just by $P_i$ when there is no ambiguity. Observe also that $z-A$ denotes the matrix $ zI_n-A\,$ \textcolor{black}{where $I_n \in \mathbb{C}^{n\times n}$ represents the identity matrix}. It is worth noting that when $A$ is Hermitian, $P_\lambda$ represents an orthogonal projection onto the eigenspace associated with $\lambda$. The Riesz projector has many interesting properties and it can be defined in more general Banach spaces, see \cite[chapter 11, pp 418-425]{riesz2012functional}\cite[pp 273--275]{kozlov1999differential}\cite[chapter 2, pp 63-64]{jeribi2015spectral} and references therein.\\

Another important integral used in this document is the Cauchy integral formulation for matrix functions.
Suppose that $f(z)$ is an analytic function inside and on a closed contour $\Gamma$ which encloses the entire spectrum of A. Then we have, see \cite[chapter 6,pp 427,Theorem 6.2.28]{horn1994topics}\cite[chapter 9,528,529]{golub2013matrix}\cite[chapter 1, pp 2-4]{higham2008functions},
\begin{align}
    f(A) = \displaystyle \frac{1}{2\pi i} \oint\limits_{\Gamma} f(z) (z-A)^{-1} \, \mathrm{d}z\,.
    \label{cauchy}
\end{align}
The complex matrix $(z-A)^{-1}$, if it exists, is known as the resolvent of $A$ and denoted as $R(z)=R(z,A)$. 

The \textit{adjugate} matrix  of $A$, which is denoted by  $\mbox{Adj}A$, is the complex matrix having entries as $(\mbox{Adj}A)_{ij} = (-1)^{i+j} M_{ji}(A)$, where $M_{ji}(A)$ is the $(j,i)$ minor of $A$. Here, the minor means the determinant of the $(n-1)\times (n-1)$ submatrix obtained by deleting the $j$-th row and the $i$-th column from matrix A. The adjugate of a matrix satisfies the following relation (\cite[Theorem 2, pp 119]{grossman1994elementary})
\begin{align}
   A  \,(\mbox{Adj}A) = \mbox{det}(A) I_n\,.
    \label{adjugateRelation}
\end{align}
% Finding the adjugate from its definition could be computationally expensive because it involves the computation of $n^2$ determinants of order $n-1$. In \cite{stewart1998adjugate} it is discussed an algorithm to compute the adjugate of a matrix $A$ using \eqref{adjugateRelation}, which can be applied even in some cases when $A^{-1}$ is ill-conditioned.\\ 

Let us introduce the characteristic polynomial of $A$ as follows:
\begin{equation}\label{eq:charpol}
    p_A(z) = p(z) =  \mbox{det}(z-A) = z^n+\alpha_1z^{n-1}+\dots+\alpha_{n-1}z+\alpha_{n}.
 \end{equation}
It follows from equation \eqref{adjugateRelation} that if the complex matrix $z-A$ is invertible,
\begin{equation}
    B(z)\equiv \mbox{Adj}(z-A) = (z-A)^{-1} p(z)\,.
    \label{eq:Adjz}
\end{equation}
Throughout the rest of this text, the notation $B(z)$ will be used to abbreviate $\mbox{Adj}(z-A)$.

It is well known that  every square complex matrix $A$ admits a Jordan decomposition. When $A$ is diagonalizable, the Jordan decomposition reduces to diagonalization. The Jordan canonical form is usually introduced in the following way. Let $A$ be a complex square matrix of order $n$ and $\{\lambda_i\}_{1 \leq i \leq m}$ represent the $m$ distinct eigenvalues of $A$, each with an algebraic multiplicity of $n_i$, $1 \leq i \leq m$. The characteristic polynomial $p(z)$ can be expressed as 
\[
p(z)=\prod_{i=1}^m (z-\lambda_i)^{n_i}.
\]
Recall that if $h$ is a positive integer, we denote by $I_{h} \in \mathbb{C}^{h\times h}$   the identity matrix of size $h \times h$. The root subspaces $\{W_i\}_{1\leq i\leq m}$ corresponding to each eigenvalue $\lambda_i$ can be defined as $W_{\lambda_i}=W_i = \mbox{Ker}(A-\lambda_i)^{n_i}$, along with the corresponding projections $P_i$ onto the spaces $W_i$. Then there exists a matrix $V$ such that $A = V J V^{-1}$, where $J$ is a block diagonal matrix, see \cite[Chapter 2]{gohberg2006invariant}\cite{grossman1994elementary,horn1994topics},
\begin{align*}
J=\operatorname{Diag}\Big(J\left(\lambda_{1}\right), \ldots, J\left(\lambda_{m}\right)\Big), 
\end{align*}
where for $i=1, \ldots, m$, the matrix $J\left(\lambda_{i}\right) \in \mathbb{C}^{n_{i} \times n_{i}}$ is given by
$$
J\left(\lambda_{i}\right)=\operatorname{Diag}\left(J_{1}\left(\lambda_{i}\right), \ldots, J_{g_{i}}\left(\lambda_{i}\right)\right), 
$$
with $ J_{k}\left(\lambda_{i}\right)=\lambda_{i} I_{n_{i,k}}+\widetilde{N}_{i k} \in \mathbb{C}^{n_{i,k} \times n_{i,k}}$ and 
$$
\widetilde{N}_{i k}=\left[\begin{array}{cc}0 & I_{n_{i,k}-1} \\ 0 & 0\end{array}\right]$$
for $k=1, \ldots, g_i$. The value $g_i$ is known as the geometric multiplicity of the eigenvalue $\lambda_i$. Numbers $\{n_{i,k}\}_{k=1,..,g_i}$ are the partial multiplicities of the corresponding eigenvalue $\lambda_i$. Moreover, 
\begin{equation}\label{eq:g_and_sik}
\sum_{k=1}^{g_i} n_{i,k} = n_i,    
\end{equation}
which is the algebraic multiplicity of $\lambda_i$. Furthermore, each $J_k(\lambda_i)$ is related to a Jordan chain.\\

Define $\widetilde{N}_{i}=\operatorname{Diag}\left(\widetilde{N}_{i 1}, \ldots, \widetilde{N}_{i g_{i}}\right) \in \mathbb{C}^{n_{i} \times n_{i}}$, put $U=V^{-1}$ and partition $U$ and $V$ according to the sizes of the diagonal blocks in $J$; that is,
$$
V=\left[\begin{array}{llll}
V_{1} & V_{2} & \cdots & V_{m}
\end{array}\right]\quad \mbox{ and } \quad U=\left[\begin{array}{c}
U_{1} \\
U_{2} \\
\vdots \\
U_{m}
\end{array}\right],
$$
where $V_{i} \in \mathbb{C}^{n \times n_{i}}$ and $U \in \mathbb{C}^{n_{i} \times n}$. The columns of $V_{i}$ are generalized eigenvectors of $\lambda_i$ and span the corresponding root subspace $W_i$.

Finally, for $i=1, \ldots, m$, introduce
$$
P_{i}=V_{i} U_{i} \quad \mbox{ and } \quad N_{i}=V_{i} \widetilde{N}_{i} U_{i}.
$$
Note that $U_{i} V_{i}=I_{n_{i}}$ for $i=1, \ldots, m$. 
%The following properties follow:
%\begin{itemize}
   % \item 
    Also that, for $i=1, \ldots, m$ we have $ P_{i}^{2}=P_{i}$, $P_{i} N_{i}=N_{i} P_{i}=N_{i}$ and $\mbox{rank}(P_i) = n_i$.
%    \item We can write $A=\displaystyle \sum_{i=1}^{m}\left(\lambda_{i} I_{n_{i}}+N_{i}\right) P_{i}.$
   % \item It holds $ \displaystyle %\sum_{t=1}^{m} P_{t}=I_{n}$ and also $P_{i} P_{j}=\delta_{i j} P_{i}$ for $
    %i, j=$ $1, \ldots, m$.
%\end{itemize}
We have the following result
(see \cite{gohberg2006invariant,grossman1994elementary,horn1994topics}).
\begin{theorem}\label{thm:TheoJordan}
Let $A$ be a complex square matrix of order $n$ and its $m$ distinct eigenvalues $\{\lambda_i\}_{1\leq i\leq m}$ with corresponding algebraic multiplicities $\{n_i\}_{1\leq i\leq m}$. Additionally, let $\{W_i\}_{1\leq i\leq m}$ be the root subspaces and $\{P_i\}_{1\leq i\leq m}$ the corresponding projections. Then $\mathbb{C}^n = \bigoplus W_i$ and there exist $N_i$, nilpotent operators of degree $d_i \leq n_i$, $1\leq i\leq m$,  such that $P_i N_i 
 = N_i P_i$ for all $i$ (we write the commutator $[N_i,P_i]=0$, for short), and 
\begin{equation}    \label{TheoJordan}
A = \sum_{i=1}^m \left( N_i + \lambda_i \right) P_i
\end{equation}
and the projectors satisfy $ \sum_{i=1}^m P_i = I_n$ and 
\begin{align} 
    P_i P_j &= \delta_{ij} P_i\quad \text{for all $1\leq i,j\leq m$}.
\end{align}
\end{theorem}

Note that if $A$ is a diagonalizable matrix (e.g.,  Hermitian matrices), we then have $N_i=0$ for $1\leq i\leq m$.\\

Coming back to the adjugate matrix, in \cite{1parisse:hal-00003444}, B. Pairisse and M.Vaughan express 
\begin{equation}\label{PV1aroundzero}
B(z) = z^{n-1}+z^{n-2}C_1+\dots+zC_{n-2}+C_{n-1}    
\end{equation}
and find a method  based on the Faddev Algorithm in order to compute $C_\ell$, $1\leq \ell\leq n-1$, using the recurrence relation
\begin{equation}\label{eq:recurrece}
C_0=I_n, \quad C_\ell = AC_{\ell-1} -\frac{1}{\ell}\mbox{tr}(AC_{\ell-1})I_n\quad \mbox{ for } \ell=1,2,\dots,n-1.
\end{equation}
Note that $
   \frac{d^{k} B(z)}{dz^{k}}|_{0} =k!C_{n-k-1}.$
A byproduct of this recurrence relation is the coefficients of the characteristic polynomial of $A$ in \eqref{eq:charpol}, that is, 
\begin{equation}\label{eq:recurrecechar}
    \alpha_{\ell}=-\frac{1}{\ell}\mbox{tr}(AC_{\ell-1}), \quad \ell=1,2,\dots,n-1.
\end{equation}
They also write the expansion centered around a given eigenvalue, say $\lambda_i$, 
\begin{equation}\label{PV1}
B(z) = \sum_{0\leq k \leq n-1} B_k(\lambda_i) (z-\lambda_i)^k.    
\end{equation}
Note that 
 $$
   \frac{d^{n_i-1} B(z)}{dz^{n_i-1}}\Bigr|_{\lambda_i} =(n_i-1)!B_{n_i-1}(\lambda_i).
 $$
Concerning this expansion they show that $\mbox{Span}(B_{n_i-1})$ is the root space associated with $\lambda_i$, i.e., the subspace $W_i$ in our notation; see \cite[Theorem 2]{1parisse:hal-00003444}. Additionally, they propose an algorithm to compute the Jordan canonical  form that we summarize here as follows. Given a matrix $A$, proceed as follows:
\begin{enumerate}
    \item Compute the matrices $C_1,C_2,\dots,C_{n-1}$ and the coefficients of the characteristic polynomial $\alpha_1,\alpha_2,\dots,\alpha_{m}$ using the Faddev Algorithm for matrix $A$ (see \eqref{eq:recurrece} and \eqref{eq:recurrecechar}).
    \item Compute the matrices $B_k(\lambda_i)$, $0\leq k\leq n-1$ and $1\leq i\leq n$. This could be done using the Horner division of $B(z)$ (and its derivatives) by $z-\lambda_i$. Another alternative is to compute them directly from the Faddev algorithm for $A-\lambda_i$.
    \item Compute the corresponding eigenvectors and generalized eigenvectors from the matrices $B_k(\lambda_i)$.
\end{enumerate}
This last step requires joint and careful columnwise Gauss elimination  for all the matrices $B_k(\lambda_i)$ that preserves the following structure, 
\begin{eqnarray*}
    (A-\lambda I_n) B(\lambda_i)&=&0  \\
    (A-\lambda I_n) B_1(\lambda_i)&=&B(\lambda_i) \\
    \vdots\\
 (A-\lambda I_n) B_{n_i-1}(\lambda_i)&=&B_{n_i-2}(\lambda_i)\\  
 (A-\lambda I_n) B_{n_i}(\lambda_i)-B_{n_i-1}(\lambda_i)&=&
 -\prod_{j\not= i} (\lambda_i-\lambda_j)^{n_i}I_n.\\
\end{eqnarray*}
 
In this manuscript, we present a different proof of some of these results using the Riesz projector. Furthermore, it is also found that $B_{n _i-1}$ is proportional to the projection to the root space $W_i$. Additionally, the matrix $N_i$ is expressed in terms of the derivatives of $B$, that is, it is provided a matrix representation  of the nilpotent matrix of the Theorem \ref{thm:TheoJordan}. As a related result, we mention that recently in \cite{denton2023eigenvectors} M. Franchi has shown some relations between the Riesz projection and the Jordan structure of a matrix.\\

\section{Higer order derivatives of the adjugate matrix}\label{sec:derivatives}
%%%%%%%%%%%%%%%%%%%%%%%%%%%%%

In order to prove our main results, first, based on formula \eqref{eq:Adjz}, let us express $(A-z)^{-1}$ in terms of the matrices $P_i$ and $N_i$, $1\leq i \leq m$. For this, we need the following general algebraic lemma.

\begin{lemma}
Let $\{X_i\}_i$ and  $\{Y_i\}_i$ be families of matrices such that $[X_i,Y_i] = \mathbf{0}$ and $X_i$ is invertible, for all i.

Assume that $\sum_j Y_j = I_n$ and that for all $i,j$ $Y_i Y_j = \delta_{ij}Y_i$. Then if $X=\sum X_i Y_i$  we can write $X^{-1} = \sum X_i^{-1} Y_i.$    
\label{lemma:LemmaInversa}
\end{lemma}
\begin{proof} We show this by direct computations as follows,
\begin{align*}
    X \sum_j X_j^{-1} Y_j &= \sum_i X_i Y_i  \sum_j X_j^{-1} Y_j\,, & & \\
    &= \sum_i \sum_j X_i Y_i X_j^{-1} Y_j\,, & & \\
    &= \sum_i \sum_j X_i Y_i Y_j X_j^{-1}\,, & & \text{because $[X_i,Y_i]=0$, hence $[X_i^{-1},Y_i]=0$,}\\
    &= \sum_i \sum_j X_i \delta_{ij} Y_i X_j^{-1}\,,& & \text{because $Y_i Y_j =\delta_{i,j}Y_i$,}\\
    &= \sum_i Y_i X_i \sum_j \delta_{ij} X_j^{-1}\,,& & \text{because $[X_i,Y_i]=0$,}\\
    &= \sum_i Y_i X_i X_i^{-1}\,,\\
    &= \sum_i Y_i\,,& &\\
    &= I_n\,.& & 
\end{align*}
 \end{proof}
 %%%%
 %%%
 %TEXT HABLAR SOBRE FUNIONES DE MATRICES
 %RELAICONAR CON LA FOMRULA INTEGRAL DE CACUHY DE ARRIBA
 %tipo .. NOte que si la integral se toma bajo un simple curve que enciierra todos los gamma_i, se reduce a la formula integral de cauchy y por lo tanto ....

We still need some additional results. The following result relates the complex Cauchy integral of $f(z) R(z)$ with the matrices $N_i$ and $P_i$. Recall that analytic functions evaluated on the matrix can be defined with the integral form of Cauchy \eqref{cauchy}; see \cite{higham2008functions}. Similar results are obtained when the contour $\Gamma$ does not enclose the entire spectrum, as if the function were restricted to the subspaces $W_i$ of the eigenvalues enclosed by $\Gamma$. 
We present complete proof for the sake of completeness.  
\begin{theorem}
     Let $f(z)$ be an analytic function and let $\Gamma_r \subset \mathbb{C}$ be a positively oriented rectifiable Jordan curve that does not contain eigenvalues of $A$ and encloses a region $D_r$ containing only one eigenvalue $\lambda_r$ of $A$. Then
     \begin{align*}
     \displaystyle \frac{1}{2\pi i} \oint\limits_{\Gamma_r} f(z) (z-A)^{-1} \, \mathrm{d}z =  f(\lambda_r + N_r) P_r\,.
     \end{align*}
     \label{Theo:CauchyProjection}
\end{theorem}
\begin{proof}
First note that if $z \in \mathbb{C}$ does not belong to the spectrum of $A$, then
%\sum_i z P_i + \sum_i\left( -\lambda_i - N_i\right) P_i
\begin{align*}
    z-A &= z-\sum_i \left( \lambda_i + N_i\right) P_i\, ,& & \text{By \eqref{TheoJordan},}\\
    z-A &= \sum_i z P_i + \sum_i \left( -\lambda_i - N_i\right) P_i\, ,& & \text{because $I_n= \sum_i P_i$ according to theorem \ref{thm:TheoJordan},}\\
    z-A &= \sum_i \left( (z-\lambda_i)-N_i \right) P_i\, ,& & \text{ }\\
    z-A &= \sum_i (z-\lambda_i) \left( I_n-\frac{1}{(z-\lambda_i)}N_i \right) P_i\, ,& & \text{}\\
    (z-A)^{-1} &= \sum_i \frac{1}{(z-\lambda_i)}
   \left( \sum_{l=0}^{n_i-1} \frac{N_i^l}{(z-\lambda_i)^l} \right) P_i\, ,& & \text{Applying Lemma \ref{lemma:LemmaInversa} with $X_i = (z-\lambda_i) \left( I_n-\frac{1}{(z-\lambda_i)}N_i\right)$. }
   \label{eq:inverseZmA}
\end{align*}
In the last step, since $N_i$ is nilpotente, then $I_n-\frac{1}{(z-\lambda_i)}N_i$ is \textcolor{black}{an invertible complex matrix} and its inverse is $\sum_{l=0}^{n_i-1} \frac{N_i^l}{(z-\lambda_i)^l}$, since $n_i$ is greater than or equal than the degree of $N_i$. \textcolor{black}{The last formula for $(z-A)^{-1}$ is well known. See \cite[p.~521. Eq 6.6.3]{horn1994topics}}.\\
%%%%
%%%%%%%%%%
After multiplying by $f(z)$, integrating over $\Gamma_r$ and applying the Cauchy formula, we have:
\begin{align*}
    \frac{1}{2\pi i} \oint\limits_{\Gamma_r} f(z) (z-A)^{-1} \, \mathrm{d}z &=  \sum_i \sum_{l=0}^{n_i-1} 
 \left( \frac{1}{2\pi i} \oint\limits_{\Gamma_r} \frac{f(z)}{(z-\lambda_i)^{l+1}} \, \mathrm{d}z  \right) N_i^l P_i\, ,& & & \\
  &=  \sum_i \sum_{l=0}^{n_i-1} 
 \delta_{r,i} \frac{f^{(l)}(\lambda_i)}{l!} N_i^l P_i\, ,& & \text{since if $r\neq i$, 
 $\frac{f(z)}{(z-\lambda_i)^{l+1}}$ is analytical in $D_r$,} & \\
 &=  \sum_{l=0}^{n_r-1} \frac{f^{(l)}(\lambda_r)}{l!}N_r^l P_r\, ,& & & \\
 &=  \sum_{l=0}^{\infty} \frac{f^{(l)}(\lambda_r)}{l!}N_r^l P_r\, ,& & \text{because $N_r^l = 0$ for $l\geq n_r$,} & \\
 &= f(\lambda_r + N_r) P_r\, .& &  &
\end{align*}
\end{proof}
%%%%%%%%%%%
As a remarkable fact, note that for $l>0$, $N_r^l P_r = (N_r P_r)^l$, because $[P_r,N_r]=0$ and $P_r^l = P_r$ given that $P_r$ is a projection. Then, $f(\lambda_r + N_r) P_r$ can be expressed as a function of $\lambda$,$P_r$ and $N_r P_r$.\\

  Note that with these results and  \eqref{cauchy}, we have that
\[
    f(A) = \sum_i f(\lambda_i + N_i) P_i\, ,
\]
and 
\[
   f(A)P_r = f(\lambda_r + N_r) P_r\, .
\]
\textcolor{black}{This is expected because the root subspace $\mbox{Ker}(A-\lambda_i)^{n_i}$ is $A-$invariant, $f(A)-$invariant and $f(AP_i) = f(A)P_i$, see \cite[Theorem 3.3.1]{gohberg2006invariant}.}
\textcolor{black}{Thus, Theorem \ref{Theo:CauchyProjection}\footnote{While this statement seems to be a well-known fact, the authors of this paper could not provide a precise reference of exactly this statement of the result.} says that $ \frac{1}{2\pi i} \oint\limits_{\Gamma_i} f(z) (z-A)^{-1} \, \mathrm{d}z$ is the restriction of $f(A)$ to $\mbox{Ker}(A-\lambda_i)^{n_i}$.} 

Before stating some results, let us define 
\begin{equation}
q_i(z) = p(z)/(z-\lambda_i)^{n_i} =\prod_{j\not =i }^m (z-\lambda_j)^{n_j}\,,
\label{eq:def:qi}
\end{equation}
for $1\leq i\leq m$, i.e, $q_i$ is the multiplication of the other factors of $p$ with $q_i(\lambda_i) \neq 0 $. Note that
\begin{equation}
q_i(\lambda_i) = \frac{1}{n_i!}\frac{d^{n_i}p(z)}{dz^{n_i}}\Bigr|_{\lambda_i}\,.
\label{eq:def:qi_lambdai}
\end{equation}
The previous theorem can be applied taking into account the equation \eqref{eq:Adjz}, and we can obtain the following result.

\begin{theorem} \label{thmmain} Let $B(z)$ be defined in \eqref{eq:Adjz} and
$q_i$ be defined in \eqref{eq:def:qi}. Then, for $0\leq s < n_i-1$
    \begin{equation}
    \frac{d^{s} B(z)}{dz^{s}}\Bigr|_{\lambda_i} = s! N_i^{n_i-1-s}\, \textbf{   } q_i\left(N_i + \lambda_i\right)  P_i.
        \label{The:BzeqPbefore}
    \end{equation}
\textcolor{black}{In addition, when} $s=n_i-1$, we have 
    \begin{equation}
    \frac{d^{n_i-1} B(z)}{dz^{n_i-1}}\Bigr|_{\lambda_i} = (n_i-1)! \, \textbf{   } q_i\left(N_i + \lambda_i\right)  P_i.\label{The:BzeqP}
    \end{equation}
\label{The:Bzmain}
\end{theorem}
\begin{proof}
We know that 
\begin{equation}
 B(z) = p(z) (z-A)^{-1}= (z-\lambda_i)^{n_i} q_i(z) (z-A)^{-1}  
 \label{eq:step}
\end{equation}
and therefore,
\begin{align}
    \frac{B(z)}{(z-\lambda_i)^{s+1}} =(z-\lambda_i)^{n_i-1-s} q_i(z) (z-A)^{-1}.
    \label{eq:AuxiliarProofDerivativesB}
\end{align}
Since $B(z) = \mbox{Adj}(z-A)$, each one of its entries is a cofactor of $z-A$ and hence $B(z)$ is an analytic function (see also \eqref{PV1aroundzero}). \textcolor{black}{Then, applying the Cauchy formula to the left-hand side of \eqref{eq:AuxiliarProofDerivativesB} and the Theorem \ref{Theo:CauchyProjection} to the right hand side, we get equations \eqref{The:BzeqPbefore} and \eqref{The:BzeqP}.}
\end{proof}
\textcolor{black}{Note that in the previous theorem, if $n_i = 1$, then there is no integer $s$, such that $0 \leq s < n_i-1$ and the equation \eqref{The:BzeqPbefore} cannot be applied. However, equation \eqref{The:BzeqP} is still valid. Furthermore, the reason why equation \eqref{The:BzeqP} is not a particular case of the former is that it does not make sense when $N_i = 0$ and $s = n_i-1$.}\\

We also have the following interesting observation on the rank of the derivatives of the adjugate matrix. \textcolor{black}{Since the only distinct eigenvalue of $q_i(N_i + \lambda_i)$ is $q_i(\lambda_i)\neq 0$, then $q_i(N_i + \lambda_i)$ is an invertible complex matrix. Therefore, for integers $0\leq s < n_i-1$, Theorem \ref{The:Bzmain} implies that 
$$ \mbox{rank}\left(\frac{d^{s} B(z)}{dz^{s}}\Bigr|_{\lambda_i}\right) = \displaystyle \mbox{rank}\left(N_i^{n_i-1-s}\right) = 
n_i-\mbox{dim}(\mbox{Ker}(A-\lambda_i)^{n_i-1-s})=
n_i-\sum_{\ell=1}^{n_i-1-s}{\#}\Big\{k \, :\, 1 \leq k \leq m, \quad n_{i,k} \geq \ell \Big\},$$ where $\{n_{i,k}\}_{k=1\ldots g_i}$ are the partial multiplicities corresponding to the eigenvalue $\lambda_i$ and $\#$ represents the cardinality of the set. For the last identity, we refer to \cite[Proposition 2.2.6]{gohberg2006invariant}. Additionally, it follows from equation \eqref{The:BzeqP} that $$\mbox{rank}\left(\frac{d^{n_i-1} B(z)}{dz^{n_i-1}}\Bigr|_{\lambda_i}\right) = \mbox{rank}(P_i) = n_i.$$ } We conclude that the rank of the first $n_i-1$ derivatives of the adjugate matrix does not decrease with the order of the derivatives up to $n_i$ which is the rank of $\frac{d^{n_i-1} B(z)}{dz^{n_i-1}}\Bigr|_{\lambda_i}$.\\

%\textcolor{blue}{Finally, for $s\geq n_i$, by the equation \eqref{eq:AuxiliarProofDerivativesB}, we have that in $\displaystyle \frac{d^{s} B(z)}{dz^{n_i-1}}\Bigr|_{\lambda_i}$ appears derivatives of $q_i(z)$ evaluated in $\lambda_i$ and additional terms containing the others matrices $N_j$ and $P_j$, with $j\neq i$.}\\

\textcolor{black}{Note that if the geometric multiplicity $g_i$ of $\lambda_i$ is such that $g_i \geq 2$, then the degree $d_i$ of $N_i$ satisfies $d_i \leq n_i - 1$. Therefore, $N_i^{n_i-1} = 0$. Thereby, taking $s=0$ in the previous theorem in equation \eqref{The:BzeqPbefore} we obtain the  well-known fact  (\cite{Castillo2021OnAF})}
    \begin{align}
        \mbox{Adj}(\lambda_i-A) = 0\,.
    \end{align}
\textcolor{black}{We can  generalize this result as follows. Let $\overline{n}_i=\max_{1\leq k\leq  g_i} n_{i,k}$ be the maximum partial multiplicity of $\lambda_i$ and let $\underline{n}_i=\min_{1\leq k\leq  g_i} n_{i,k}$ be the corresponding minimum partial multiplicity. Therefore, if 
%$ s\geq 0$ is such that $d_i \leq  n_i-1-s$ ( that is 
$0 \leq s \leq n_1-1-\overline{n}_i$
%), 
then $N_i^{n_i-1-s} = 0$. Therefore, by equation \eqref{The:BzeqPbefore}, $\frac{d^{s} B(z)}{dz^{s}}\Bigr|_{\lambda_i} = 0$. In addition, it holds for all $s$ such that $0\leq s \leq (g_i-2)\underline{n}_i-1$, when the corresponding geometric multiplicity $g_i \geq 2$. Since $\underline{n}_i(g_i-1) \leq n_i - \overline{n}_i $, see \eqref{eq:g_and_sik}}. We have the following results.
\begin{corollary}
    Let $A$ be a $n \times n$ square matrix, $\lambda_i$ an eigenvalue of $A$. Let $n_i$ and $g_i$ be the algebraic and geometric multiplicities of $\lambda_i$, respectively. Let $\overline{k}_i$ and $\underline{n}_i$ be the maximum and minimum partial multiplicities of $\lambda_i$, respectively. If $s$ is an integer such that $0\leq s \leq n_i-1-\overline{n}_i$, then
    \begin{align}
    \frac{d^{s} \mbox{Adj}(z-A)}{dz^{s}}\Bigr|_{\lambda_i} = 0.
    \label{eq:ZerosDerivatives}
    \end{align}
    In particular, if $ g_i\geq 2$ and $0\leq s \leq \underline{k}_i(g_i-1)-1$, then the equation \eqref{eq:ZerosDerivatives} holds.
\end{corollary}
\textcolor{black}{
Now, when the geometric and algebraic multiplicities coincide (e.g., every eigenvalue of a diagonalizable matrix), then the nilpotent matrix $N_i$ associated with $\lambda_i$ is the zero matrix, and we have the following results.}
\begin{corollary}
Let $A$ be a $n \times n$ square matrix, $B(z) = \mbox{Adj}(z-A)$, $\lambda_i$ an eigenvalue of $A$ with geometric multiplicity equal to the algebraic multiplicity $n_i$, then for every integer $0\leq s < n_i-1 $, we have
    \begin{align*}  \label{The:TMcRecover}
    \frac{d^{s} B(z)}{dz^{s}}\Bigr|_{\lambda_i} &= 0\,.
    \end{align*}
We also have (see \eqref{eq:def:qi_lambdai})
    \begin{equation}
    \frac{d^{n_i-1} B(z)}{dz^{n_i-1}}\Bigr|_{\lambda_i} = (n_i-1)! \, \textbf{   } q_i\left(\lambda_i\right)  P_i
    = \frac{1}{n_i}\frac{d^{n_i}p(z)}{dz^{n_i}}\Bigr|_{\lambda_i}\ P_i.
    \end{equation}
In particular:%, we have the following statements.
\begin{itemize}
    \item If $\lambda_i$ is simple, then 
    \begin{align}
    B(\lambda_i) =\frac{d p(z)}{dz}\Bigr|_{\lambda_i} P_i\,.
    \label{Corol:BeqPi}
     \end{align}
    \item If $A$ is hermitian and $\lambda_i$ is a simple eigenvalue, with a unit eigenvector $z_i$, then
\begin{align}
    B(\lambda_i)&=\frac{d p(z)}{dz}\Bigr|_{\lambda_i} z_i z_i^*\,.
\end{align}
\end{itemize}
\end{corollary}

As before, when $A$ is Hermitian and $\lambda_i$  is simple, then $n_i = 1$, $N_i = 0$ and $\mbox{rank}(P_i) = 1$, hence $P_i$ has a matrix representation of the operator $z_i z_i^*$, and the formula is reduced to the one of the Thompson and McEnteggert Theorems stated next.
\begin{theorem}[Thompson and McEnteggert]
Suppose that $A$ is a hermitian matrix and $\lambda_i$ \textcolor{magenta}{is} a simple eigenvalue with unit eigenvector $z_i$, then
\begin{align}
    \mbox{\normalfont Adj}(\lambda_i - A) = \frac{d p(z)}{dz}\Bigr|_{\lambda_i} z_i z_i^*\,.
\end{align}
\label{eq:TheoremTM}
\end{theorem}
For this result, we refer to \cite{iserles2002acta,parlett1998symmetric,thompson1968principal}.
See also \cite{Castillo2021OnAF} for a generalization of this identity to include any matrix with entries in an arbitrary field that is stated as follows. Let $\lambda_i$ be a simple eigenvalue of a matrix $A$ with right and left eigenvectors $v_i$ and $z_i$, respectively. Then ( \cite[Remark 2.10.]{Castillo2021OnAF}),
\begin{equation}
\mbox{\normalfont Adj}(\lambda_i - A) = \frac{d p(z)}{dz}\Bigr|_{\lambda_i}\frac{1}{z^*_iv_i}v_i z_i^*\,.
\end{equation} 
\textcolor{black}{Observe that $\frac{1}{z^*_iv_i}v_i z_i^*$ is an oblique projection on $v_i$ in the direction orthogonal to $z_i$}. 
Theorem \ref{The:Bzmain} includes this case since $n_i=1$ implies $N_i=0$ and $\mbox{rank}(P_i)=1$. In this case, $P_i$ is an oblique projection in the space generated by the associated eigenvector $v_i$.  The extension of the results to a more general algebraic setting (as in \cite{Castillo2021OnAF})  will be presented elsewhere.\\

For the following discussion, assume that $\lambda_i$ has geometric multiplicity $g_i=1$ and $n_i > 1$.
 Then there is only one Jordan block associated with $\lambda_i$ and only one,  up to scaling, right (and also one left)  eigenvector. Select arbitrarily $v_i$ and $u^*_i$ right and left eigenvectors, respectively. 
We have $N_i^{n_i} = 0$ and $N_i^{n_1-1} \neq 0$. Thus, $N_i^{n_i-1} q_i(N_i + \lambda_i) = q_i(\lambda_i) N_i^{n_i-1}$ (since 
$z^{n_i-1}q_i(z+\lambda_i)=q_i(\lambda_i)
z^{n_i-1}+q'(\lambda_i)z^{n_i}+\dots $). Therefore, by equation \eqref{The:BzeqPbefore}
\begin{equation}
    B({\lambda_i}) = q_i(\lambda_i) N_i^{n_i-1}\, P_i = q_i(\lambda_i) N_i^{n_i-1}.
\end{equation}
But $N_i^{n_i-1}$ is a rank 1 matrix 
that can be written as $N_i^{n_i-1} = \alpha v_i u_i^*$ with $\alpha \in \mathbb{C}$ since the $\mbox{Image}(N_i^{n_i-1})=\mbox{span}\{v_i\}$. Let $x_i$ be a generating vector of $v_i$ and and $y_i$ be a generating vector of $u_i$, that is, 
$$(\lambda_i - A)^{n_i-1} x_i = (-1)^{n_i-1} N_i^{n_i-1} x_i = v_i 
\quad \mbox{ and }\quad  y_i^*(\lambda_i - A)^{n_i-1} = u_i^*.$$ 
We have $u_i^*x_i \neq 0$ and 
$\alpha=\frac{(-1)^{n_i-1}}{u_i^*x_i}$ since $v_i = (-1)^{n_i-1} N_i^{n_i-1} x_i = (-1)^{n_i-1} \alpha (u_i^* x_i)v_i$.  Therefore $$N_i^{n_i-1} = (-1)^{n_i-1} \frac{v_i u_i^*}{u_i^*x} = (-1)^{n_i-1} \frac{v_i u_i^*}{y^* (\lambda_i - A)^{n_i-1} x}.$$ Thus we have recovered, for the particular case of complex matrices, the result stated in \cite[C. 2.8]{Castillo2021OnAF}. Finally observe that  $ \alpha^2(u_i^*v_i)v_iu_i^*=  (\alpha v_iu_i^*)(\alpha v_iu_i^*)   =N_i^{n_i-1} N_i^{n_i-1}=0$, thus, $u_i^*v_i = 0$. In summary, we have the following result. See Example 
\ref{example3} for an illustrative example.

\begin{corollary}
    Let $A\in \mathbb{C}^{n\times n}$. Suppose that $\lambda_i$ is an eigenvalue of $A$ with geometric multiplicity $g_i=1$ and algebraic multiplicity $n_i >  1$. Let $v_i,u_i$ be corresponding right and left eigenvectors of $A$, with eigenvalues $\lambda_i$ and $\lambda_i^*$. In addition, let $x_i,y_i$ be generating right and left eigenvectors of $v_i$ and $u_i$. respectively. Then
    \begin{equation}
        \mbox{Adj}(\lambda_i-A) = (-1)^{n_i-1} q_i(\lambda) \frac{v_i u_i^*}{y^* (\lambda_i - A)^{n_i-1} x}\,.
    \end{equation}
    Moreover, the left and right vectors are orthogonal, that is, $v_i^*u_i = 0$.
    \label{Cor:Mul1lAlg}
\end{corollary}

% $q_i(N_i+\lambda_i)P_i$ can be considered as a operator form $W_i$ onto $W_i$, because
% \begin{align*}
%     q_i(N_i+\lambda_i)P_i = q_i(A)P_i = \Pi_{j\neq i} (A- \lambda_j)^i P_i = \Sigma a_r (A-\lambda_i)^r
% \end{align*}
% with $a_0 \neq 0$ because $\lambda_i$ is not a root of $q_i$. Since $W_i$ is invariant under $(A-\lambda_i)$, then so under $\Sigma a_r (A-\lambda_i)^r = q_i(N_i+\lambda_i)$. Now, it is invertible (as an operator form $W_i$ to $W_i$), because by the Jordan theorem, $q_(A) = V q(J) V^{-1}$, since $J$ is an upper triangular matrix, then $q(J)$ is an upper triangular matrix where its diagonal entries are $q(J)_{t,t} = q(J_{t,t})$. When projector $P_i$ acts on $q_(A)$, then just remains the Jordan blocks associated with $\lambda_i$, therefore $q_(A)P_i=P_iq_(A)P_i$ (seen as an operator from $W_i$ to $W_i$) is an upper triangular matrix, where its diagonal entries all are $q(\lambda_i)\neq 0$, hence, it is invertible. 
% 
% Now, from \ref{The:BzeqP}, it is deduced that $span(B^{(n_i-1)}(\lambda_i)) \subset W_i$, because $B^{(n_i-1)}(\lambda_i)$ is the multiplication of $P_i$ and other matrix. Now, if $y \in W_i$, then, since $q_i(A_i)=q_i(N_i+\lambda_i)P_i$is invertible as an operator from $W_i$ to $W_i$, then, it does exist $x\in W_i$ such that $q_i(A_i)x = y$, then $B^{(n_i-1}(\lambda_i) x = q_i(A_i)x = y$, hence $W_i \subset span(B^{(n_i-1}(\lambda_i))$ and therefore $W_i = span(B^{(n_i-1}(\lambda_i))$, thus we recover the result of \cite{1parisse:hal-00003444}.

To finish this section we express $P_i$ and $N_i$ in terms of the derivatives of $B(z)$ and $q_i(z)$. Considering $f(z) = 1$ and applying Theorem \ref{Theo:CauchyProjection}, we have
\begin{align*}
    P_i=\frac{1}{2\pi i} \oint\limits_{\Gamma_i} (z-A)^{-1} \, \mathrm{d}z = \frac{1}{2\pi i} \oint\limits_{\Gamma_i} \frac{B(z)}{(z-\lambda_i)^{n_i} q_i(z)} \, =\frac{1}{(n_i-1)!}\frac{d^{n_i-1} }{dz^{n_i-1}} \left( \frac{B(z)}{q_i(z)}\right)\Bigr|_{\lambda_i}\,.
\end{align*}
If $n_i = 1$, then $N_i = 0$. However, if $n_i\geq 2$, taking $f(z) = z-\lambda_i$ and applying Theorem \ref{Theo:CauchyProjection}, then,
\begin{align*}
    N_i = \frac{1}{2\pi i} \oint\limits_{\Gamma_i} (z-\lambda_i) (z-A)^{-1} \, \mathrm{d}z &= \frac{1}{2\pi i} \oint\limits_{\Gamma_i} \frac{B(z)}{(z-\lambda_i)^{n_i-1} q_i(z)}=\frac{1}{(n_i-2)!}\frac{d^{n_i-2} }{dz^{n_i-2}} \left( \frac{B(z)}{q_i(z)}\right)\Bigr|_{\lambda_i}\,.
\end{align*}
Hence we have the following result.

\begin{theorem} \label{thmmain} Let $B(z)$ be defined in \eqref{eq:Adjz} and
$q_i$ be defined in \eqref{eq:def:qi}. Then
\begin{equation}
    P_i=\frac{1}{(n_i-1)!}\frac{d^{n_i-1} }{dz^{n_i-1}} \left( \frac{B(z)}{q_i(z)}\right)\Bigr|_{\lambda_i}\,.
    \label{The:BzQzeqP}
\end{equation}
Moreover, if $n_i\geq2$, then
\begin{align}
    N_i &= \frac{1}{(n_i-2)!}\frac{d^{n_i-2} }{dz^{n_i-2}} \left( \frac{B(z)}{q_i(z)}\right)\Bigr|_{\lambda_i}\,.
    \label{The:BzeqNPeq}
\end{align}
\end{theorem}
Therefore, one can recover the matrices $P_i$ and $N_i$ from the derivatives of cofactors terms. Therefore we  can compute functions $f(A)$ with finite polynomials because of the fact that $N_i$ is nilpotent.

\section{Illustrative examples}\label{sec:examples}
In this section, we present particular illustrative examples to help the reader to fix ideas related to the results presented in the previous section. The code can be found in 
\url{https://github.com/jrubianom/Jordan_From_Derivatives}.
\subsection{Example 1}
Let A be the following hermitian matrix
\begin{align}
    A &=  
\begin{pmatrix}
  \begin{matrix}
    1 & -1 & 1\\
    -1 & 1 & -1 \\
    1 & -1 & 1
  \end{matrix}
\end{pmatrix}\,.
\end{align}
The characteristic polynomial is $p(z) = z^2 (z-3)$. Then the eigenspace related with $\lambda = 0$  has dimension 2 and the eigenspace related with $\lambda = 3$ has dimension 1. The adjugate matrix is 
$$B(z)=\left(\begin{matrix}z^{2} - 2 z & - z & z\\- z & z^{2} - 2 z & - z\\z & - z & z^{2} - 2 z\end{matrix}\right)$$ and 
its  derivative is written as
$$
B^{\prime}(z)=
\left(\begin{matrix}2 z - 2 & -1 & 1\\-1 & 2 z - 2 & -1\\1 & -1 & 2 z - 2\end{matrix}\right).
$$
Evaluating $B^{\prime}(z)$ at $z=0$ we get
$$
    P_0 =  
    \frac{1}{-3}B'(0)=
\begin{pmatrix}
    \frac{2}{3} & \frac{1}{3} & -\frac{1}{3}\\
    \frac{1}{3} & \frac{2}{3} & \frac{1}{3}\\
    -\frac{1}{3} & \frac{1}{3} & \frac{2}{3}
\end{pmatrix}
$$
and  evaluating $B(z)$ at $z = 3$
$$
P_3 = 
\frac{1}{3^2 \cdot 0!}B(3)=
\begin{pmatrix}
  \begin{matrix}
    \frac{1}{3} & -\frac{1}{3} & \frac{1}{3}\\
    -\frac{1}{3} & \frac{1}{3} & -\frac{1}{3}\\
    \frac{1}{3} & -\frac{1}{3} & \frac{1}{3}
  \end{matrix}
\end{pmatrix}\,.
$$
One can check the following relations
\[    P_0 + P_3 = I_n\,, \quad  P_0 P_3 = P_3 P_0 = 0\,,  \quad   P_0^2 = P_0\,, \quad 
    P_3^2 = P_3\,,  \quad AP_0 = 0\,   \quad    AP_3 = 3 P_3    \quad \mbox{ and }  A = 3 P_3 + 0 P_0\,.
\]
These relations are basically the spectral theorem for this specific matrix. Note that $\mbox{Tr}(P_0) = 2$ and $\mbox{Tr}(P_3) = 1$, i.e, the dimensions of the eigenspaces, as it should be (the trace of a projection is the dimension of its range). See \cite{Daniel, galantai2013projectors}.

As a remark related to our results, we observe that the Jordan canonical form for this matrix is given by 
$$
    V =  
\begin{pmatrix}
  \begin{matrix}
    1 & -1 & 1\\
    1 & 0 & -1\\
    0 & 1 &1
  \end{matrix}
\end{pmatrix}
\quad \mbox{ and }  \quad
   J     =
\begin{pmatrix}
  \begin{matrix}
    0 & 0 & 0\\
    0 & 0 & 0\\
    0 & 0 &3
  \end{matrix}
\end{pmatrix}
$$
We verify with the eigenvectors corresponding to the first eigenvalues, $V_0= V(:,[1,2])$ that,
$P_0=V_0*(V_0^TV_0)^{-1}V_0^T$
and with $V_3=V(:,3)$ that 
$
P_3=V_3*(V_3^TV_3)^{-1}V_3^T.
$

\subsection{Example 2}
Consider now the following matrix
$$
   A 
    =
\begin{pmatrix}
  \begin{matrix}
    0 & 1 & 0 & 0\\\\11 & 6 & -4 & -4\\\\22 & 15 & -8 & -9\\\\-3 & -2 & 1 & 2
  \end{matrix}
\end{pmatrix}.
$$
The characteristic polynomial is $p(z) = \left(z - 1\right)^{2} \left(z + 1\right)^{2}$. The root space related to $\lambda = 1$ has dimension 2, the same for $\lambda = -1$. The adjugate matrix is given by
$$
   B(z)
    =
\begin{pmatrix}
  \begin{matrix}
    z^{3} + 9 z - 10 & z^{2} + 6 z - 7 & 4 - 4 z & 4 - 4 z\\\\11 z^{2} - 10 z - 1 & z^{3} + 6 z^{2} - 7 z & - 4 z^{2} + 4 z & - 4 z^{2} + 4 z\\\\22 z^{2} + 16 z - 26 & 15 z^{2} + 10 z - 17 & z^{3} - 8 z^{2} - 7 z + 10 & - 9 z^{2} - 6 z + 11\\\\- 3 z^{2} - 6 z - 3 & - 2 z^{2} - 4 z - 2 & z^{2} + 2 z + 1 & z^{3} + 2 z^{2} + z
  \end{matrix}.
\end{pmatrix}
$$
Since the algebraic multiplicity is $2$ in both cases, we compute the first derivative of $B(z/q_1(z)$ that is given by, 

$$
   \left(\frac{B(z)}{q_1(z)}\right)^{\prime}
    =
\begin{pmatrix}
  \begin{matrix}
    \frac{- 2 z^{3} - 18 z + 3 \left(z + 1\right) \left(z^{2} + 3\right) + 20}{\left(z + 1\right)^{3}} & \frac{4 \left(5 - z\right)}{z^{3} + 3 z^{2} + 3 z + 1} & \frac{4 \left(z - 3\right)}{\left(z + 1\right)^{3}} & \frac{4 \left(z - 3\right)}{\left(z + 1\right)^{3}}\\\frac{8 \left(4 z - 1\right)}{z^{3} + 3 z^{2} + 3 z + 1} & \frac{z^{3} + 3 z^{2} + 19 z - 7}{z^{3} + 3 z^{2} + 3 z + 1} & \frac{4 \left(1 - 3 z\right)}{z^{3} + 3 z^{2} + 3 z + 1} & \frac{4 \left(1 - 3 z\right)}{z^{3} + 3 z^{2} + 3 z + 1}\\\frac{4 \left(7 z + 17\right)}{z^{3} + 3 z^{2} + 3 z + 1} & \frac{4 \left(5 z + 11\right)}{z^{3} + 3 z^{2} + 3 z + 1} & \frac{z^{3} + 3 z^{2} - 9 z - 27}{z^{3} + 3 z^{2} + 3 z + 1} & - \frac{12 z + 28}{z^{3} + 3 z^{2} + 3 z + 1}\\0 & 0 & 0 & 1
  \end{matrix}.
\end{pmatrix}
$$
Evaluating at $z=1$, we have
$$
   P_1 = \left(\frac{B(z)}{q_1(z)}\right)^{\prime}\Bigr|_{z=1}
    =
\begin{pmatrix}
  \begin{matrix}
    3 & 2 & -1 & -1\\\\3 & 2 & -1 & -1\\\\12 & 8 & -4 & -5\\\\0 & 0 & 0 & 1
  \end{matrix}
\end{pmatrix}.
$$
%%%%%%
%%%%%%
%%%%%%%55
%%%%%%%%%%%55
%%%%%%%%%%%%%%%%
%%%%%%%%%%%
%%%%%%%%%%
%
and 
$$
   \left(\frac{B(z)}{q_{-1}(z)}\right)^{\prime}
    =
\begin{pmatrix}
  \begin{matrix}
    \frac{z^{2} - 2 z - 11}{z^{2} - 2 z + 1} & - \frac{8}{\left(z - 1\right)^{2}} & \frac{4}{\left(z - 1\right)^{2}} & \frac{4}{\left(z - 1\right)^{2}}\\- \frac{12}{\left(z - 1\right)^{2}} & \frac{z^{2} - 2 z - 7}{z^{2} - 2 z + 1} & \frac{4}{\left(z - 1\right)^{2}} & \frac{4}{\left(z - 1\right)^{2}}\\\frac{12 \left(3 - 5 z\right)}{z^{3} - 3 z^{2} + 3 z - 1} & \frac{8 \left(3 - 5 z\right)}{z^{3} - 3 z^{2} + 3 z - 1} & \frac{z^{3} - 3 z^{2} + 23 z - 13}{z^{3} - 3 z^{2} + 3 z - 1} & \frac{8 \left(3 z - 2\right)}{z^{3} - 3 z^{2} + 3 z - 1}\\\frac{12 \left(z + 1\right)}{z^{3} - 3 z^{2} + 3 z - 1} & \frac{8 \left(z + 1\right)}{z^{3} - 3 z^{2} + 3 z - 1} & - \frac{4 z + 4}{z^{3} - 3 z^{2} + 3 z - 1} & \frac{z^{3} - 3 z^{2} - 5 z - 1}{z^{3} - 3 z^{2} + 3 z - 1}
  \end{matrix}
\end{pmatrix}.
$$
Evaluating at $z=-1$, we have
$$
   P_{-1} = \left(\frac{B(z)}{q_{-1}(z)}\right)^{\prime}\Bigr|_{z=-1}
    =
\begin{pmatrix}
  \begin{matrix}
    -2 & -2 & 1 & 1\\\\-3 & -1 & 1 & 1\\\\-12 & -8 & 5 & 5\\\\0 & 0 & 0 & 0
  \end{matrix}
\end{pmatrix}.
$$

%%% ME FALTA AÑADIR LOS N_i P_i

The Jordan canonical form for this matrix is given by 
$$
    V =  
\begin{pmatrix}
  \begin{matrix}
    \frac{1}{3} & - \frac{2}{3} & 0 & \frac{1}{4}\\- \frac{1}{3} & 1 & 0 & \frac{1}{4}\\\frac{1}{3} & 0 & \frac{1}{4} & 1\\0 & 0 & - \frac{1}{4} & 0
  \end{matrix}
\end{pmatrix}
\quad \mbox{ and } \quad 
   J 
    =
\begin{pmatrix}
  \begin{matrix}
    -1 & 1 & 0 & 0\\\\0 & -1 & 0 & 0\\\\0 & 0 & 1 & 1\\\\0 & 0 & 0 & 1
  \end{matrix}
\end{pmatrix}.
$$
Thus,
$$
   \mathcal{I}_1
    =
\begin{pmatrix}
  \begin{matrix}
    0 & 0 & 0 & 0\\\\0 & 0 & 0 & 0\\\\0 & 0 & 1 & 0\\\\0 & 0 & 0 & 1
  \end{matrix}
\end{pmatrix}
\quad 
\mbox{ and }
\quad    \mathcal{I}_{-1}    =
\begin{pmatrix}
  \begin{matrix}
    1 & 0 & 0 & 0\\\\0 & 1 & 0 & 0\\\\0 & 0 & 0 & 0\\\\0 & 0 & 0 & 0
  \end{matrix}
\end{pmatrix}
$$

We verify $
    P_1 = V \mathcal{I}_1 V^{-1}$ and $
    P_{-1} = V \mathcal{I}_{-1} V^{-1}.$

Analogously, we obtain
$N_iP_i =  \left( \frac{B(z)}{q_i(z)}\right)\Bigr|_{\lambda_i}$, that is, 

$$
   N_1
    =
\begin{pmatrix}
  \begin{matrix}
    0 & 0 & 0 & 0\\\\0 & 0 & 0 & 0\\\\3 & 2 & -1 & -1\\\\-3 & -2 & 1 & 1
  \end{matrix}
\end{pmatrix}
\quad 
\mbox{ and }
\quad 
   N_{-1}
    =
\begin{pmatrix}
  \begin{matrix}
    -5 & -3 & 2 & 2\\\\5 & 3 & -2 & -2\\\\-5 & -3 & 2 & 2\\\\0 & 0 & 0 & 0
  \end{matrix}
\end{pmatrix}.
$$
We have verified the following relations, 
\begin{align*}
    &P_{1} + P_{-1} = I_n\,,\quad 
    P_{1} P_{-1} = P_{-1} P_1 = 0\,,\quad 
    P_{1}^2 = P_{1}\,,\quad 
    P_{-1}^2 = P_{-1}\,,\quad
    AP_1 = (N_1+1)P_1\,,\\ &AP_{-1} = (N_{-1}-1) P_{-1}\,,\quad 
    (N_1 P_1)^2 = (N_{-1}P_{-1})^2 = 0, \quad \mbox{ and }
    A = (N_1+1)P_1 + (N_{-1}-1) P_{-1}\, .
\end{align*}

As an example of a function of $A$,
$$\mbox{exp}(A) = \mbox{exp}(N_1+1)P_1 + \mbox{exp}(N_{-1}-1)P_{-1} = \mbox{exp}(1) P_1 + \mbox{exp}^\prime(1) N_1 P_1 + \mbox{exp}(-1) P_{-1} + \mbox{exp}^\prime(-1) N_{-1} P_{-1},$$
We verify that this coincides with $\mbox{exp}(A)$
$$
   \mbox{exp}(A)
    =
\begin{pmatrix}
  \begin{matrix}
    - \frac{7}{e} + 3 e & - \frac{5}{e} + 2 e & - e + \frac{3}{e} & - e + \frac{3}{e}\\\frac{2}{e} + 3 e & 4 \cosh{\left(1 \right)} & - 2 \cosh{\left(1 \right)} & - 2 \cosh{\left(1 \right)}\\- \frac{17}{e} + 15 e & - \frac{11}{e} + 10 e & - 5 e + \frac{7}{e} & - 6 e + \frac{7}{e}\\- 3 e & - 2 e & e & 2 e
  \end{matrix}
\end{pmatrix}.
$$

\subsection{Example 3}\label{example3}
Consider the matrix
$$
A = V_1 \begin{pmatrix}
\begin{matrix}1 & 1 & 0 & 0 & 0\\0 & 1 & 1 & 0 & 0\\0 & 0 & 1 & 0 & 0\\0 & 0 & 0 & 2 & 0\\0 & 0 & 0 & 0 & 5 i\end{matrix}
\end{pmatrix} V_1^{-1}\,,\, 
V_1 = 
\begin{pmatrix}
    \begin{matrix}
    1 & 0 & 1 & 0 & \frac{1}{2}\\
    1 & 1 & 1 & 0 & \frac{1}{2}\\
    1 & 0 & 2 & 0 & \frac{1}{2}\\
    1 & 0 & 1 & 1 & \frac{1}{2}\\
    1 & 0 & 1 & 0 & 1
    \end{matrix}
\end{pmatrix}\,.
$$
Note that the geometric multiplicity $g_1$ of the eigenvalue $1$ is $g_1=1$. Whereas the algebraic multiplicity is $n_1 = 3$.\\
The Jordan decomposition of $A^{*}$ is
$$
A^* = V_2 
\begin{pmatrix}
\begin{matrix}
1 & 1 & 0 & 0 & 0\\
0 & 1 & 1 & 0 & 0\\
0 & 0 & 1 & 0 & 0\\
0 & 0 & 0 & 2 & 0\\
0 & 0 & 0 & 0 & - 5 i
\end{matrix}
\end{pmatrix} V_2^{-1}\,,\, 
V_2 = 
\begin{pmatrix}
    \begin{matrix}1 & 1 & -2 & -1 & -1\\0 & -1 & 0 & 0 & 0\\-1 & 0 & 0 & 0 & 0\\0 & 0 & 0 & 1 & 0\\0 & 0 & 1 & 0 & 1\end{matrix}
\end{pmatrix}\,.
$$
According to Corollary \ref{Cor:Mul1lAlg}, the left and right eigenvectors of $A$ are orthogonal. Thus
$$
\left(\begin{matrix}1 & 0 & -1 & 0 & 0\end{matrix}\right)^* \left(\begin{matrix}1 & 1 & 1 & 1 & 1\end{matrix}\right)^t = 0\,.
$$
\section{Final Comments}\label{sec:conclusions}
In this manuscript, we obtained explicit formulas relating    higher order derivatives of  the adjugate matrix $\mbox{Adj}(z-A)$ to  the Jordan decomposition of the matrix $A$.  See Theorem \ref{thmmain}.
To obtain these identities we used  the Riesz projector and some results from functional calculus. 
The results presented here can be considered a generalization of the Thompson and McEnteggert theorem that relates the adjugate matrix with the orthogonal projection on the eigenspace of simple eigenvalues for symmetric matrices. See \cite{iserles2002acta,parlett1998symmetric,thompson1968principal}.
It can also be regarded as a generalized form of some identities in \cite{Castillo2021OnAF}.
They can also be viewed as a complement to some  previous results by B. Parisse, M. Vaughan in  \cite{1parisse:hal-00003444} that related derivatives of the adjugate matrix with the invariant subspaces associated with an eigenvalue. Additionally, the formulas can be regarded as general eigenvector-eigenvalue identity, see \cite{denton2022eigenvectors}.
Although this method for obtaining the eigenvectors from derivatives of the cofactors of $A-z$ and  the nilpotent matrices of the Jordan decomposition (in the base generated by the eigenvectors) may not be efficient from a numerical point of view, it may be useful for symbolic computations or  theoretical purposes.
Further generalization of the main results of this paper over other fields than complex may be carried out in the future. Applications of the techniques developed here to other problems in linear algebra related to functional calculus and invariant subspaces are the subject of current research.

\section*{Authorship contribution statement}
J.R-M. wrote an initial paper draft 
with the results and proofs presented here. J.G. Helped write the final version of the manuscript and to establish connections with existing previous results. 
\section*{Acknowledgements}
This work was developed  in the 2022-II Numerical Analysis and Finite Element Seminar at the Universidad Nacional de Colombia - Bogot\'a.
The authors thank Professor Marcus Sarkis for proposing to study the article \cite{frommer2001algebraic} that led them to discuss the spectrum of nonsymmetric matrices. 

\printbibliography
\end{document}